\documentclass[12pt,a4paper,reqno]{amsart}
\usepackage{amsmath,amsfonts,amsthm,amssymb,verbatim,enumerate,graphicx}

\usepackage{hyperref}

\usepackage[marginratio=1:1,totalwidth=15.75cm,totalheight=22.275cm]{geometry}

\newcommand*{\defeq}{\mathrel{\vcenter{\baselineskip0.5ex \lineskiplimit0pt
                     \hbox{\scriptsize.}\hbox{\scriptsize.}}}%
                     =}

\newcommand{\R}{\mathbb{R}}
\newcommand{\C}{\mathbb{C}}
\newcommand{\Z}{\mathbb{Z}}
\newcommand{\eps}{\varepsilon}

\newtheorem{theorem}{Theorem}

\begin{document}

\title{Points of convergence - music meets mathematics}
\author{Lasse Rempe}
\address{Department of Mathematical Sciences\\ University of Liverpool\\ Liverpool L69 7ZL\\ United Kingdom\textsc{\newline \indent \href{https://orcid.org/0000-0001-8032-8580}{\includegraphics[width=1em,height=1em]{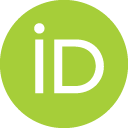} {\normalfont https://orcid.org/0000-0001-8032-8580}}}} 
\email{lrempe@liverpool.ac.uk}

\date{This is a preprint of the chapter ``Points of convergence - music meets mathematics,'' to appear in ``More UK Success Stories in Industrial Mathematics,'' ed. Philip J.\ Aston.}

\begin{abstract}\emph{Phase-locking} is a fundamental phenomenon in which coupled 
or periodically forced oscillators synchronise. The 
 \emph{Arnold family} of circle maps,
 which describes a forced oscillator, is the simplest
  mathematical model of phase-locking  and 
 has been studied intensively since its introduction in the 1960s. 
 The family exhibits regions of parameter space
  where phase-locking phenomena can be observed. A long-standing question asked whether 
  \emph{hyperbolic} parameters~-- those whose behaviour is dominated by periodic attractors, and which are therefore stable under perturbation~-- are
  dense within the family. A positive answer was given in 2015 by van Strien and the author, which implies that, no matter
  how chaotic a map within the family may behave, there are always systems with stable behaviour nearby. 
  This research was a focal point of a pioneering collaboration with 
  composer Emily Howard, commencing with Howard's residency in Liverpool's mathematics department in 2015. 
  The collaboration generated impacts on creativity, culture and society, including several 
musical works by Howard,
  and lasting influence on artistic practice through a first-of-its-kind centre for science and music. We describe the research and the collaboration, and reflect on the 
  factors that contributed to the latter's success.
  \end{abstract}

\maketitle

\section*{Introduction}

In the 17th century, the Dutch scientist 
  Christiaan Huygens discovered that two pendulum clocks,
  coupled by being mounted on the same wooden beam,
  synchronised their movements. This phenomenon,
  described by Huygens as ``a miraculous sympathy,'' 
   is called
  \emph{phase locking} (also mode locking, or entrainment): 
  interacting periodic oscillators tend to synchronise their movements
  (with the same period or with periods that are related by an integer
  multiple). Phase locking is near-ubiquitious in physical oscillators: Examples
  include the fact that the revolution period and orbital period of the moon are identical 
  and the synchronisation of fireflies. A similar effect
  occurs when one oscillator is \emph{forced}
  by another, for example in bowed musical instruments: When 
  a violin string is plucked by hand, it exhibits non-harmonic overtones; when 
  bowed, all of these overtones are forced onto the harmonic series. 
  
  A goal of ``pure'' (i.e., theoretical) mathematics is to investigate interesting
   phenomena in their most fundamental settings, to understand the 
   underlying principles. In a 1961 paper~\cite[\S~12]{arnold}, Vladimir Arnold introduced what may be
   the simplest model for phase-locking behaviour: a family of self-maps of the
   circle, motivated by the movement (in discrete time-steps) of a forced
   periodic oscillator. Known as the \emph{Arnold family} of circle maps, it can be described by the formula
    \begin{equation}\label{eqn:arnold} f_{\alpha,b}(\theta) \defeq \theta + \alpha + 
              b \sin(\theta) \quad (\operatorname{mod} 2\pi).\end{equation}
   The angle $\theta$ represents
   a point on the circle $S^1 = \R/2\pi \Z$. 
   The number $\alpha$ is a rotation parameter (also an element of $S^1$); if $b=0$, then $f_{a,b}$ is simply the rotation by
   an angle of $\alpha$. 
    The final term in~\eqref{eqn:arnold} is a forcing term, which must also be 
     a function of $\theta$, and thus periodic. We use the
     simplest possible periodic forcing term, namely a multiple of $\sin$.\footnote{%
  Arnold used $\cos$ instead of $\sin$, which is equivalent.} The parameter
     $b\in [0,\infty)$ determines the strength of the forcing. 

  Let us fix parameters $\alpha$ and $b$ and write $f = f_{\alpha,b}$. Given a
   starting state $\theta_0\in\R/2\pi\Z$, we think of $\theta_1 = f(\theta_0)$ as
    determining the state of our dynamical system after one time step.  
   So $\theta_2 = f(\theta_1)$ is the state after two steps, and after $n$ steps: 
   \[ \theta_{n} = f(\theta_{n-1}) = f(f(... f(\theta_0)\dots)) \quad\text{($n$ times)}.
      \]
   The sequence $(\theta_n)_{n=0}^{\infty}$ 
   is called the \emph{orbit} of $\theta_0$ under $f$. 
   If $\theta_n = \theta_0$ for some $n>0$, then the (finite) orbit is called 
   a \emph{periodic cycle}; this cycle is \emph{stable} (under perturbation of $\theta_0$) 
   if the orbits of all nearby
   starting values converge to it.

\begin{figure}
 \begin{center}
 \includegraphics[width=\textwidth]{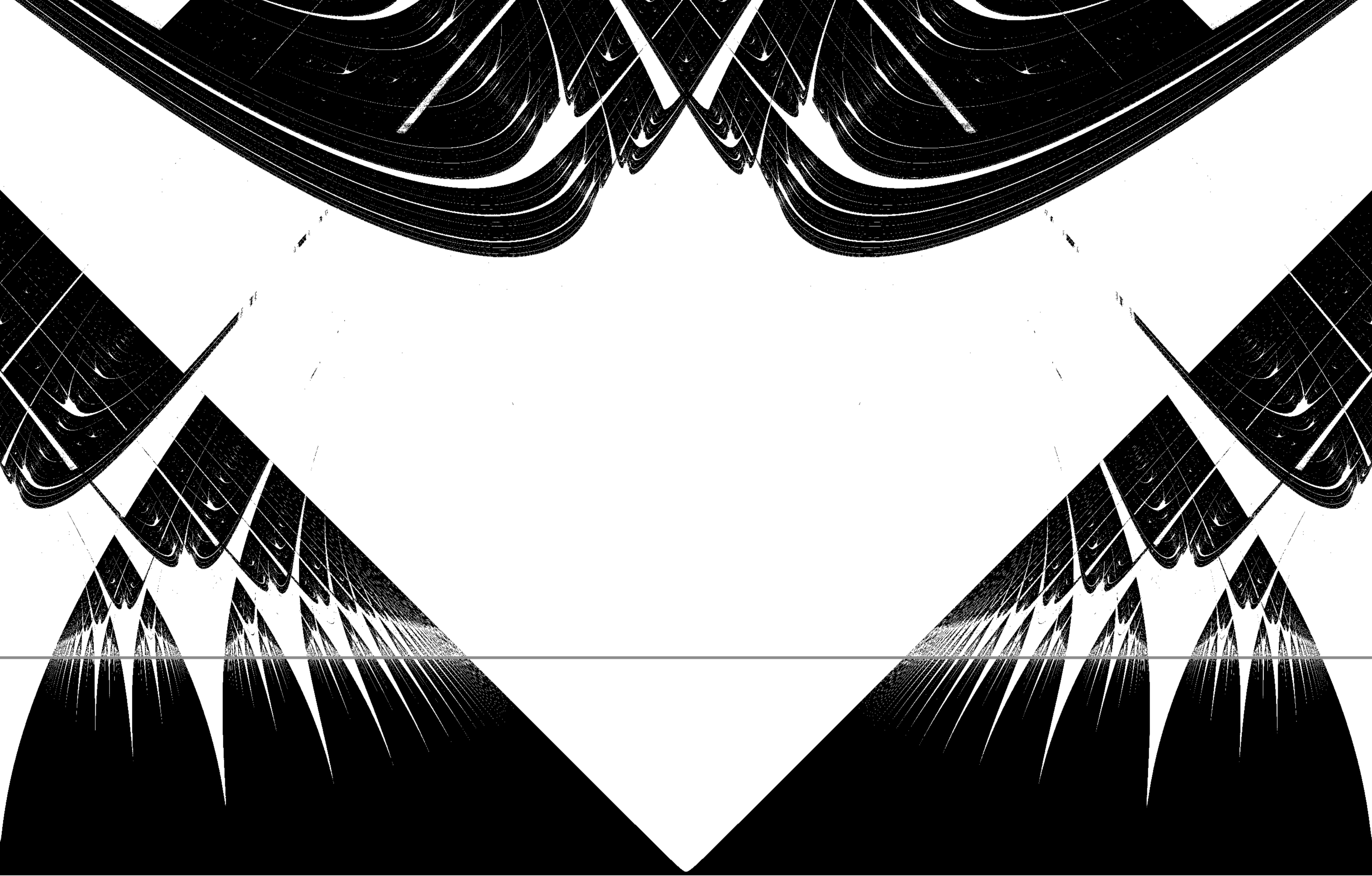}
\end{center}
\caption{Parameter space of the Arnold family. The horizontal and vertical 
 directions correspond, respectively, to the parameters $\alpha$, ranging from $-\pi$ to $\pi$ and
 $b$, ranging from $0$ to $4$. The critical line $b=1$ is shown in grey; white 
regions correspond to stable (hyperbolic) maps. Theorem~\ref{thm:main} states that these regions are dense in parameter space; the apparent presence of solid black regions arises from the
 fact that stable regions may be very thin.\label{fig:arnold}}
\end{figure}
     
     For $b<1$ (the bottom quarter of Figure~\ref{fig:arnold}), 
       the map $f_{\alpha,b}$ is a diffeomorphism of the circle. It indeed 
       exhibits phase-locking phenomena, which are well-understood: 
      For parameters in certain regions, called ``Arnold tongues''
       (the white regions), 
      all but finitely many orbits tend to a stable 
      periodic cycle. Moreover, these maps are stable under perturbations of $f$ within the 
       Arnold family: 
       a small
       change of the parameter leads to a small change in the overall behaviour of orbits,
       so $f_{\alpha,b}$ is ``phase-locked'' to the periodic orbit. 

       
    For $b>1$, the map $f_{\alpha,b}$ is no longer invertible and has two distinct
     critical points on $S^1$. The Arnold tongues begin to
     intersect, and the periodic orbit for a given tongue may bifurcate and become
     unstable. This can lead to \emph{chaotic} behaviour, where an arbitrarily small
     change of the starting state $\theta_0$ could lead to completely 
     different long-term behaviour. 
     It is thus natural to ask whether for $b>1$, there is still a dense set of parameters 
     whose behaviour is stable under perturbations of the parameter. 

   More precisely, $f=f_{\alpha,b}$ with $b>1$ is called \emph{hyperbolic} if the
    orbits of both 
    critical points
     tend to stable periodic cycles. Such $f$ is stable under 
      perturbations of the parameter~-- 
      all nearby maps are also hyperbolic~-- and almost every orbit
      converges to a stable cycle. Density of hyperbolic
      maps in parameter spaces is a central question of one-dimensional dynamics. For 
      real polyomials, it was established by Kozlovski, Shen and 
      van Strien~\cite{KSS} in 2007, answering part (b) of Smale's 11th problem. However, 
      this does not resolve the question in
      the Arnold family: a key issue is that $f$,
      when extended to the complex plane, is transcendental rather than algebraic (see below). This problem was overcome
       by van Strien and the author~\cite{density}:\footnote{%
 In fact, the results of~\cite{density} are more general, covering many
   families of transcendental entire functions and circle maps.}

   \begin{theorem}[Density of hyperbolicity]\label{thm:main}
    Hyperbolic maps are dense in the Arnold family.
     That is, given any $\alpha\in S^1$ and $b>1$, and every $\eps>0$, there 
       exist perturbed parameters $\alpha'\in S^1$ and $b'>1$ with 
       $\lvert \alpha - \alpha'\rvert,\lvert b - b'\rvert < \eps$
       such that $f_{\alpha',b'}$ is hyperbolic.
   \end{theorem}

\begin{figure}
 \begin{center}
  \includegraphics[width=\textwidth]{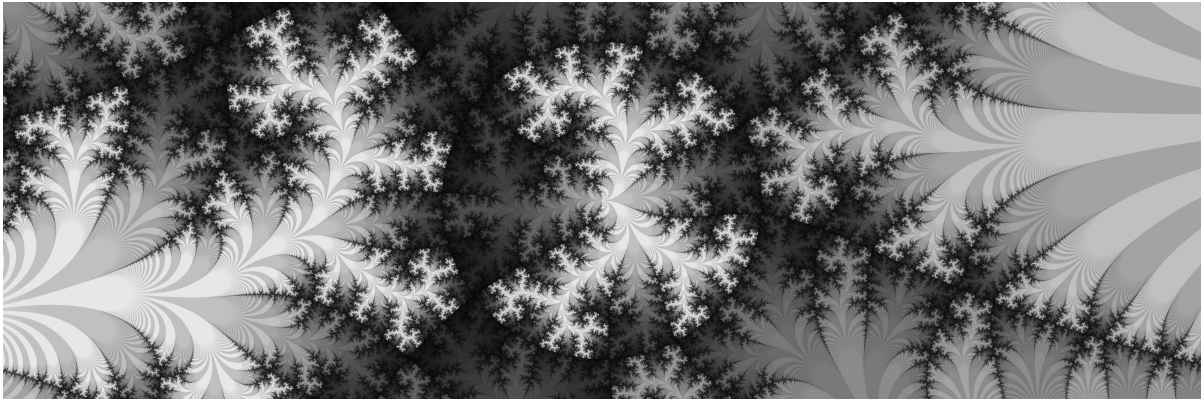}
\end{center}
\caption{The dynamics of $g_{\alpha,b}$ in the complex plane, for one choice of 
  parameters. As usual, the horizontal and vertical directions correspond to 
  the real and imaginary parts of the starting value. 
  The map is chaotic on the unit circle
  (visible in the central part of the figure), and almost every orbit is escaping. 
   Different shades of grey distinguish 
    different
    patterns with which orbits tend to the essential singularities
    at $0$ and $\infty$, highlighting the intricate dynamics of the complex extension.\label{fig:arnold-complex}}
\end{figure}

\section{Ideas in the proof of the theorem}
 To prove Theorem~\ref{thm:main}, we consider the 
  \emph{complex extension} of $f_{\alpha,b}$, allowing the state 
  $\theta$ in~\eqref{eqn:arnold} to be a complex number. Applying the change of variable 
  $z = \exp(i\theta)$, we obtain a self-map 
   of the punctured plane $\C^*= \C\setminus \{0\}$ 
   (see Figure~\ref{fig:arnold-complex}):
   \[ g_{\alpha,b}(z) = \exp( i(\theta + \alpha + b\sin(\theta))) =
     e^{i\alpha} \cdot z \cdot \exp\left(\frac{1}{2}\bigl( z - \frac{1}{z}\bigr)\right).  \] 

 This function $g_{\alpha,b}$ 
     has isolated singularities at $0$ and $\infty$, which are \emph{essential}
    (neither removable nor poles). Consequently, the behaviour
    near these points is very 
    complicated; for example, the preimages of every $z\notin \{0,\infty\}$ 
   accumulate both 
    at $0$ and at $\infty$. There is a set of starting values of positive area 
   which are \emph{escaping}; i.e., their orbits
    accumulate only on the essential singularities $0$ and $\infty$.

    A crucial step in the proof is to establish
    \emph{rigidity}: a map in the Arnold family cannot be deformed,
    by a small change of parameters, to one with 
    the same qualitative dynamical behaviour, 
    except through certain well-understood mechanisms. In our
    setting, a new difficulty arises: We must
    exclude the existence of deformations arising from the set of
    escaping points. This issue does not arise for 
    polynomials, as treated in~\cite{KSS}, which
   have no essential singularities. It  
    is overcome by using 
    techniques developed for studying the behaviour of transcendental entire functions
    near $\infty$ \cite{boettcher}. 

To prove the theorem, we then begin
   with a map $f_{\alpha,b}$ that is not hyperbolic; recall that this map 
     has two critical points and two critical values.  
    It follows from the above-mentioned
   rigidity statement that we 
     may perturb $(\alpha,b)$ slightly to a parameter $(\alpha_1,b_1)$ for which the
    orbit of one 
    critical value passes through a
    critical point. The set of parameters satisfying this critical relation 
    forms an analytic variety of dimension $1$. We can make another perturbation,
   \emph{within this variety}, to create a second critical relation. (Otherwise, maps within
    this variety would have the same qualitative dynamical behaviour as 
    $f_{\alpha_1,b_1}$, which contradicts rigidity.) 

  Thus we obtain parameters $(\alpha',b')$, arbitrarily close to 
    $(\alpha,b)$, for which both critical values are eventually mapped
    to a critical point. This means that each critical point is mapped to a periodic cycle
    containing a critical point; such a cycle is 
    necessarily stable. Thus $f$ is hyperbolic, and the proof of the theorem is complete. 

\section{A musical collaboration} 

 In 2015, Liverpool's mathematics department 
  hosted award-winning composer Emily Howard as a Leverhulme 
  Artist in Residence. Howard previously used mathematical and scientific ideas in her
  compositions, but the mathematical discussions during the
  residency   led to the use of frontier mathematical research, rather than 
   classical mathematical principles, in her creative process for the first time.

\begin{figure}[t]
 \begin{center}
  \includegraphics[width=\textwidth,height=1.3cm]{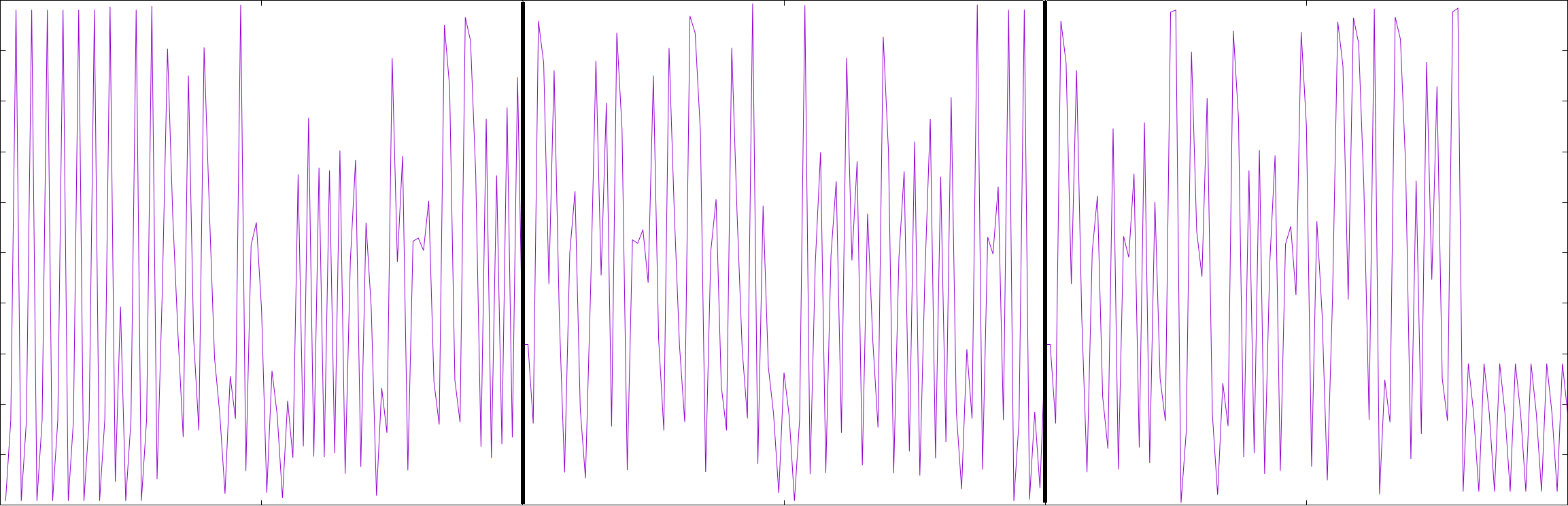}\\
  \includegraphics[width=\textwidth,height=1.3cm]{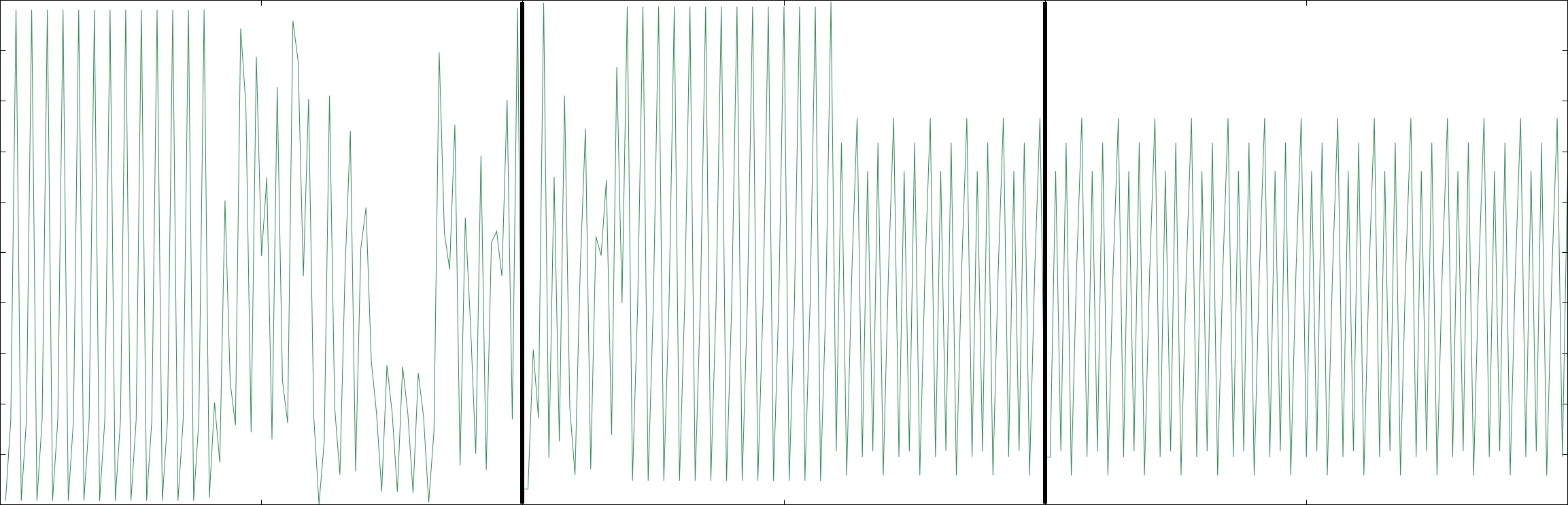}
\end{center}
 \caption{\label{fig:orbit}The data set used in Howard's composition 
   \emph{Orbit 2a} first follows two orbits of an unstable function in the Arnold family,
     starting near an unstable periodic cycle, but quickly separating
    due to the chaotic dynamics. As in the proof of Theorem~\ref{thm:main}, 
    the function is slightly perturbed once, changing the parameters by less than $10^{-5}$, to create an attracting cycle. A second, even smaller, perturbation yields 
     another attracting cycle, and thus a hyperbolic function. (A member 
    of the Arnold family with two
    attracting cycles is necessarily hyperbolic.)} 
\end{figure}

 A particular focus of discussions during the residency 
  were the article~\cite{density} and Theorem~\ref{thm:main}. Howard 
  challenged the 
  dynamics group to 
  create sets of numeric data, encapsulating key ideas of the work.
  Two datasets were created by the author and Alexandre Dezotti
   (see Figure~\ref{fig:orbit}); the resulting 
  composition, \emph{Orbits} \cite{orbits}, is a direct creative
  response to the 
  data. Other pieces
   (\emph{Leviathan}, \emph{Threnos} and \emph{Chaos or Chess}) written
  during this period also used the discussions around~\cite{density} as pivotal
  creative input, establishing enduring principles for Howard's approach to composition.   

  For example, \emph{Torus} was a BBC Proms co-commission with the Royal Liverpool Philharmonic Orchestra, 
     first performed in 2016 at the Royal Albert Hall to a sell-out audience. It was inspired by Howard's discussions
      with mathematicians in Liverpool, including 
     those with the dynamical systems group, as well as conversations with
      Liverpool geometer  
      Anna Pratoussevitch. The ideas of chaotic motion, and of a small perturbation
       that changes the fundamental nature of a system, formed an important part of Howard's compositional approach
       for this piece. \emph{Torus} was hailed by the Guardian as ‘one of this year’s finest new works’, and recognised
       with a 2017 British Composer Award. Two further geometry-inspired orchestral
      pieces followed: \emph{sphere} (2017) was commissioned by the Bamberg Symphony 
     Orchestra and broadcast on BBC Radio 3, and 
   \emph{Antisphere} was commissioned to open the London
     Symphony Orchestra's 2019/20 season at the Barbican under
     Sir Simon Rattle. 

   In 2017, the Royal Northern College of Music established a new 
     dedicated centre for Practice and Research in Science and Music (PRiSM), with
    Howard as Director. PRiSM is building on the 
    approach developed during the Liverpool residency, bringing together 
    scientists, composers and performers for mutual benefit. 


    Due to the specialised nature of mathematical research,
    it can be challenging to communicate non-superficial ideas and concepts 
    to those outside the specific area of research, let alone outside of mathematics. There 
    are two factors  
    that contributed strongly to successfully creating 
    a dialogue between music and research mathematics. We shall discuss them briefly
   here, as they may be instructive for collaborative projects of a similar nature.

The first is the presence of a shared language: Howard is an Oxford graduate in
 mathematics and computing, while the author is an amateur orchestral musician. Howard's 
  familiarity with mathematical terminology allowed for detailed discussions of 
  mathematical
  ideas, which were reflected in the resulting compositions. For instance,
  in the data underlying \emph{Orbits}, there appear repeating patterns 
  arising from both stable and unstable cycles. Though they 
   are indistinguishable from the data
  alone, Howard chose to 
  represent their differing nature musically as a result of the 
  mathematical discussions. Likewise, the perturbations of parameters
  are imperceptible in the data; knowing about their significance,
  Howard marked them as recognisable musical events. It appears
  likely that, for a collaboration like this to succeed, a common language,
   if not already present, needs to be carefully developed. 

A second key factor in the success of the collaboration 
   was the careful choice of a specific piece of research as the basis of
   discussions. There are several reasons why Theorem~\ref{thm:main} 
   provided fertile ground for developing new creative thinking: 
\begin{enumerate}
  \item It lends itself to expression in simple and general terms: no matter how chaotic the 
      system, there is always stability nearby. 
  \item Phase-locking is known and relevant to musicians, e.g. through
       the synchronisation of 
       linked metronomes or the elimination of non-harmonic overtones. 
   \item Aspects of the proof are philosophically intriguing: the complex plane~-- invisible in the formulation of the problem~-- plays a crucial part in the proof. 
   \item The systems in question can be used to design data series as input into and inspiration for the creative process.
\end{enumerate}

\subsection*{Acknowledgements.} The research in~\cite{density} was supported
  by Rempe's EPSRC Grants EP/E017886/1 and EP/E052851/1, Rempe's
 Philip Leverhulme Prize
 and van Strien's ERC Advanced Grant. Howard's residency was funded by
  a Leverhulme Artist in Residence award.

\end{document}